\documentstyle[12pt]{article}
\newcommand{\twist}[3]{{#1}${\,\scriptscriptstyle {#3}}\atop\raise9pt
           \hbox{$\scriptstyle\oo$} ${#2}}
\newcommand{\ve}{\varepsilon}
\newcommand{\be}{\begin{eqnarray}}
\newcommand{\ee}{\end{eqnarray}}
\newcommand{\n}{\nonumber }
\newcommand{\oo}{\otimes}
\newcommand{\bt}{\beta}
\newcommand{\la}{\lambda}
\newcommand{\si}{\sigma}

\newcommand{\al}{\alpha}

\begin{document}

\begin{titlepage}
\begin{center}
{\Large \bf Universal $R$-matrix for esoteric quantum group}
\end{center}
\vspace{0.5cm}
\begin{center}
 {\bf P. P. Kulish}
\end{center}
\vspace{0.5cm}
\begin{center}
St.Petersburg Department of the Steklov
Mathematical Institute,\\
Fontanka 27, St.Petersburg, 191011,
Russia;\\ ( kulish@pdmi.ras.ru )
\end{center}
\vspace{0.5cm}
\begin{center}
 {\bf A. I. Mudrov}
\end{center}
\vspace{0.5cm}
\begin{center}
Department of Theoretical Physics,
Institute of Physics, St.Petersburg State University,
St.Petersburg, 198904, Russia\\(aimudrov@dg2062.spb.edu)
\end{center}
\vspace{1cm}

\begin{abstract}
The universal $R$-matrix for a class of esoteric (non-standard)
quantum groups ${\cal U}_q(gl(2N+1))$ is constructed as a twisting
of the universal $R$-matrix ${\cal R}_S$ of the Drinfeld-Jimbo quantum
algebras. The main part of the twisting element
${\cal F}$ is chosen to be the canonical element of
appropriate pair of separated Hopf subalgebras (quantized
Borel's ${\cal B}(N) \subset {\cal U}_q(gl(2N+1))$), providing the
factorization property of ${\cal F}$. As a result,
the esoteric quantum group generators can be
expressed in terms of the Drinfeld-Jimbo ones.
\end{abstract}
\end{titlepage}

\section{Introduction}

Quasitriangular Hopf algebras (quantizations of universal enveloping
Lie algebras ${\cal U}(g)$) and quantum groups (deformations of
algebra
of functions on Lie groups $Fun(G)$) are the subject of active
research during the last decade originated in famous Drinfeld's
report \cite{D1}. There are different starting points of the quantum
group
theory: generators and defining relations \cite{D1,J}, $R$-matrix
or solution to the Yang-Baxter equation (FRT approach) \cite{RTF},
deformation quantization or star product \cite{FSt}.
Although up to now there is no complete transformation theory
of quantum groups,  particular studies of twistings are of great
importance \cite{F&C}.

In this paper  a twisting element ${\cal F}$ is constructed giving the
esoteric quantum group of \cite{FG1,FG2} as a deformation  of
the Drinfeld-Jimbo quantum group (the standard one)
and the corresponding universal $R$-matrix from that of quantum
${\cal U}_q(gl(2N+1))$. The question of  relation between ${\cal
R}_{FG}$ and ${\cal R}_{S}$ via twisting was already discussed in
Refs.
\cite{H9506018,H9609029,JC9702028} within the FRT-approach and in
terms of
matrix 2-cocycles $\chi$ on quantum groups, which are the images of
"universal" twisting elements in the fundamental representation of
the
quantum algebra: $\chi(T^i_j,T^k_l)=F^{ik}_{jl}$.
Below we are dealing with
a quasitriangular Hopf algebra ${\cal A}$ only fixing appropriate
Hopf subalgebras in it. The latter ones can be considered as
mutually dual and the twisting cocycle ${\cal F}$ is given by the
corresponding canonical element $\sum e_i\oo e^i \in {\cal B}(N)\oo
{\cal B}(N)^*_{op}$ (the subscript "op" means the opposite
multiplication).

The Letter is organized as follows. After reminding briefly the basic
material on twisting of Hopf algebras (Sec.2), we construct the universal
twist for the ${\cal U}_q(gl(3))$ case, when the Cremmer-Gervais $R$-matrix
(found in \cite{CG} for $gl(N)$) coincides with the esoteric one.
It is shown that knowing the universal twist one can express the
FRT-approach generators in terms of the original quantum algebra
generators. Sec.4 is devoted to esoteric ${\cal U}_q(gl(2N+1))$
for general $N$. The Letter is concluded by outlining few possible
applications of twisting element.

\section{Twisting of Hopf algebras}
\label{Sec2}

A Hopf algebra ${\cal A}(m,\Delta,\ve,S)$ with multiplication
$m\colon {\cal A}\oo{\cal A}\to {\cal A}$, coproduct
$\Delta\colon {\cal A}\to {\cal A}\oo{\cal A}$, counit
$\ve\colon {\cal A}\to C$, and antipode
$S\colon {\cal A}\to {\cal A}$ (see definitions in
Refs.\cite{D1,RTF,ChP})
can be transformed \cite{D2} with an invertible element
${\cal F}\in {\cal A}\oo {\cal A}$, ${\cal F} =
\sum f_i^{(1)}\oo f_i^{(2)}$ into a twisted one
${\cal A}_t(m,\Delta_t,\ve,S_t)$. This Hopf algebra ${\cal A}_t$
has the same multiplication and
counit maps but the twisted coproduct and antipode
$$
\Delta_t(a)={\cal F}\Delta_t(a){\cal F}^{-1},\quad
S_t(a)=vS_t(a)v^{-1},\quad v=\sum f_i^{(1)}S(f_i^{(2)}),\quad
a \in {\cal A}.
$$
Sometimes it appears to be useful to combine twist with a
homomorphism of ${\cal A}$. The twisting element has to satisfy  the
identities
\be
(\ve \oo id)({\cal F}) = (id \oo \ve )({\cal F})=1,
\ee
\be
{\cal F}_{12}(\Delta \oo id)({\cal F}) =
{\cal F}_{23}(id \oo \Delta)({\cal F}),
\label{TE}
\ee
where the first one is just
a normalizing condition and follows from the second relation
modulo a non-zero scalar factor.

A quasitriangular Hopf algebra  ${\cal A}(m,\Delta,\ve,S,{\cal R})$
has additionally an element ${\cal R}\in {\cal A}\oo{\cal A}$
(a universal $R$-matrix) satisfying \cite{D1}
\be
(\Delta \oo id)({\cal R})={\cal R}_{13}{\cal R}_{23}, \quad
(id \oo \Delta)({\cal R})={\cal R}_{13}{\cal R}_{12}.
\label{UR}
\ee
The coproduct $\Delta$ and its opposite $\Delta^{op}$
are related by the similarity transformation (twisting)
with ${\cal R}$
$$\Delta^{op}(a)={\cal R}\Delta(a){\cal R}^{-1},\quad  a \in {\cal
A}.$$
A twisted quasitriangular quantum algebra
${\cal A}_t(m,\Delta_t,\ve,S_t,{\cal R}_t)$  has the
twisted universal $R$-matrix
\be
{\cal R}_t=\si({\cal F})\,{\cal R}\,{\cal F}^{-1},
\label{Rt}
\ee
where $\si$  means permutation of the tensor factors:  $\si (f\otimes
g) = (g\otimes f)$.

Although, in principle, the possibility to quantize an arbitrary Lie
bialgebra has been proved  \cite{EK}, an explicit formulation of
Hopf operations remains a nontrivial task. In particular, the
knowledge of explicit form of the twisted cocycle is a rare case even
for classical universal enveloping algebras, despite of advanced
Drinfeld's theory \cite{D3}. Most of such explicitly known twisting
elements have the factorization property with respect to
comultiplication (cf.(\ref{UR}))
$$(\Delta \oo id) ({\cal F}) = {\cal F}_{23}{\cal F}_{13}
\quad \mbox{or}
\quad (\Delta \oo id) ({\cal F}) = {\cal F}_{13}{\cal F}_{23}
$$
and similar involving  $(id \oo \Delta)$.
To satisfy the twist equation,
these identities are combined with additional requirements
${\cal F}_{12}{\cal F}_{23}={\cal F}_{23}{\cal F}_{12}$
or the Yang-Baxter equation on ${\cal F}$ \cite{RSTS,R}.

In particular, a twisting element can be used to construct a
nontrivial tensor product of Hopf algebras ${\cal A}$  and ${\cal B}$
\cite{RSTS,M}.
Given  an element ${\cal F}\in {\cal A}\oo {\cal B}$, ${\cal F}
=
\sum a_i\oo b_i$ such that
\be
(\Delta_{\cal A}\oo id) ({\cal F}) ={\cal F}_{23} {\cal F}_{13}, \quad
  (id\oo \Delta_{\cal B}) ({\cal F}) ={\cal F}_{12} {\cal F}_{13},
  \label{FP}
\ee
one can  define the twisted tensor product
\twist{${\cal A}$}{${\cal B}$}{{\cal F}} coinciding with
${\cal A}\oo{\cal B}$ as an algebra and endowed with a new
coproduct
$$
\Delta_t(a\oo b)={\cal F}_{14}(id \oo \si \oo id)
(\Delta_{\cal A}(a)\oo\Delta_{\cal B}(b)){\cal F}_{14}^{-1},
$$
where ${\cal F}_{14}=\sum a_i \oo 1 \oo 1\oo b_i\in
({\cal A}\oo{\cal B})\oo ({\cal A}\oo{\cal B})$
and the antipode
$$S(a\oo b)={\cal F}^{-1} S_{\cal A}(a)\oo S_{\cal B}(b) {\cal F}.$$
Taking  ${\cal B}={\cal A}^*_{op}$ (the dual with the opposite
multiplication) and the canonical element ${\cal F}=\sum e_i\oo e^i$
as a
twisting cocycle, one gets \cite{RSTS} the dual to the Drinfeld
quantum double
$$D({\cal A}^*)= \mbox{(\twist{${\cal A}$}{${\cal A}^*_{op}$}{{\cal
F}})}^*.$$
Let us note that the element  ${\cal F}$ can be replaced by
 $(id\oo \varphi)({\cal F})$ where  $\varphi$ is a Hopf automorphism
 of ${\cal B}$. Such a modification may be nontrivial if
${\cal A}\oo{\cal B}$ is embedded into a larger Hopf  algebra ${\cal H}$
and $\varphi$ is not extended to an automorphism of entire  ${\cal H}$.
A twisting element of this kind with appropriate Hopf
algebras ${\cal A}\,, {\cal A}^*_{op} \subset {\cal U}_q(sl(2N+1)) $
will be constructed in the next sections.

\section{The Cremmer-Gervais universal $R$-matrix for $gl(3)$}

>From the study of the quantum $sl(N)$ Toda field theory \cite{CG}
the Cremmer-Gervais solution to the Yang-Baxter equation
$R_{CG}$ was obtained, which was different from the standard one
\be
R_S&=&q\sum_i e_{ii}\oo e_{ii} +\sum_{i\not= j} e_{ii}\oo e_{jj}
                           +\omega\sum_{i< j} e_{ij}\oo e_{ji},
                           \quad q=e^\gamma,\quad \omega=q-q^{-1}.
\label{RSt}
\ee
The $R$-matrix $R_{CG}$ for $N=3$ coincides with
the Fronsdal-Galindo $R$-matrix of esoteric quantum groups
\cite{FG1,FG2}
\be
R_{CG}&=& R_{S}+(p-1)\left( e_{11}\otimes e_{22}+e_{22}\otimes
e_{33}\right) \n\\
&+&(p^{-1}-1)\left( e_{22}\otimes e_{11}+e_{33}\otimes
e_{22}\right) + ( p^2 /q - 1)e_{11}\otimes e_{33}\n\\
& +& ( q/p^2 -1)e_{33}\otimes
e_{11}+q \nu \,(e_{32}\otimes e_{12}- p^2 /q^2 e_{12}\otimes e_{32} )
\,.
\label{RCG}
\ee
Here $p$ and $\nu$ are two additional independent deformation
parameters.

The quasitriangular Hopf algebra ${\cal U}_q(sl(3))$ can be defined
by the two triples of generators $\{h_i,e_i,f_i\}$, $i=1,2,$
subjected to the relations \cite{D1,J}
\be
q^{h_i} e_j = q^{a_{ij}} e_j q^{h_i}, \quad
  q^{h_i} f_j   = q^{-a_{ij}} f_j q^{h_i},\quad
      [e_i,f_j] = \delta_{ij} \frac{q^{h_i}-q^{-h_i}}{q-q^{-1}},\n\\
e_k^2 e_l - (q+q^{-1})\, e_k e_l e_k + e_l e_k^2,\quad
f_k^2 f_l - (q+q^{-1})\, f_k f_l f_k + f_l f_k^2,\quad k\not=l,
\label{triples}
\ee
where the Cartan matrix elements are
$a_{ii}=2$, $a_{ii+1}=a_{i+1 i}=-1$.
The coproduct on these generators reads
\be
\Delta(h_i)=h_i\oo 1  + 1\oo h_i,\>\>
\Delta(e_i)=e_i\oo q^{h_i} +1\oo e_i, \>\>
\Delta(f_i)=f_i\oo 1  + q^{-h_i}\oo f_i.
\label{coprod}
\ee
Having introduced elements corresponding to the composite root
$$
e_{13}=e_1 e_2 - q e_2 e_1, \quad f_{13}=f_2 f_1 - q^{-1} f_1 f_2,
$$
one gets the universal $R$-matrix in the factorized
form \cite{Ros,KR,LS,KT,T}:
\be
{\cal R}_S = q^{t_0} \exp_{q^{-2}}(\omega e_2\oo f_2)\,
                  \exp_{q^{-2}}(\omega e_{13}\oo f_{13})\,
                  \exp_{q^{-2}}(\omega e_1\oo f_1) ,
\label{R3}
\ee
where $t_0 = \sum_{ij} (a^{-1})_{ij} h_i\oo h_j$ is the canonical
element of the Cartan subalgebra ${\cal H} \oo {\cal H} $ and
the q-exponential is
\be
  \exp_q(x)&=&\sum_{n=0}^{\infty} \frac{x^n}{[n;q]!}
  =\{\prod^{\infty}_{k=0} (1-(1-q)xq^k)\}^{-1}
  ,\quad
  [n;q]!=\frac{q^n-1}{q-1}. 
\ee
The new matrix elements of $R_{CG}$ in Eq. (\ref{RCG})
correspond to contributions of
$e_1\oo f_2$ and $f_2\oo e_1$ in the fundamental representation.
Although commuting with each other,
the elements $e_1$ and $f_2$ do not generate independent
Hopf subalgebras because $h_1$ does not commute with $f_2$ nor does
$h_2$ with $e_1$. To overcome this obstacle let us
extend ${\cal U}_q(sl(3))$ with the central element
$$
C=e_{11}+e_{22}+e_{33}, \quad h_1 = e_{11} - e_{22},
\quad  h_2 = e_{22} - e_{33}\,,
$$
and perform a diagonal twist to separate the above mentioned Hopf
subalgebras
$$
{\cal F}^{(1)}=\exp(\frac{\gamma}{2}[e_{11}\wedge e_{22}+e_{11}\wedge
e_{33}+e_{22}\wedge e_{33}])\,, \quad q=e^\gamma  \, .
$$
The twisted coproducts $\Delta_t= {\cal F} \Delta {\cal F}^{-1}$
of the generators $\tilde e_1 =e_1 q^{\frac{1}{2}(e_{11}+e_{22})}$
and $\tilde f_2 =f_1 q^{-\frac{1}{2}(e_{22}+e_{33})}$
do not contain common elements:
$$
\Delta_t(\tilde e_1)=\tilde e_1\oo q^{2 e_{11}} +1\oo \tilde e_1,
\quad
\Delta_t(\tilde f_2)=\tilde f_2\oo 1  + q^{2 e_{33}}\oo\tilde f_2.
$$
The Hopf subalgebra ${\cal B}(1)_-$ generated by $\{ e_{33},\tilde f_2\}$
appears to be dual but with opposite product ${\cal B}(1)^*_{op}$
to the Hopf subalgebra
${\cal B}(1)$ spanned by  $\{ e_{11},\tilde e_1\}$. So, the corresponding
canonical element is
\be
{\cal F}^{(2)}=\exp_{q^2} (\mu \, \tilde e_1 \oo \tilde f_2)
q^{2 e_{11}\oo e_{33}}
\label{FCG}
\ee
with independent parameter $\mu$. This element
can be used for further twisting already twisted ${\cal U}_q(gl(3))$.
With extra diagonal twist depending on $(e_{11} - e_{33}) \oo C $
we get three parameter universal $R$-matrix,
reduced to (\ref{RCG}) in the fundamental representation.

Let us briefly discuss the relations among the FRT-generators of the
Drifeld-Jimbo (standard) quantum algebra and the twisted one.
Taking the first factor of ${\cal A} \oo {\cal A}$ in the fundamental
representation, we get three $3 \times 3$ matrices
$$
F_{21} = ({\rho \oo id}) \sigma({\cal F}), \quad
L^{(+)}_S = ({\rho \oo id}) {\cal R}_S, \quad
F_{12} = ({\rho \oo id}) {\cal F} \,,
$$
entries of which are expressed in terms of the standard generators
(\ref{triples}). Multiplying these matrices one gets the $L$-matrix
of the FRT-approach
$$
L_{CG}^{(+)} = F_{21} L_S^{(+)} F_{12}^{-1}
$$
entries of which are generators of the esoteric quantum algebra,
adding
the same formulas for
$$
L_{CG}^{(-)} = ({\rho \oo id}) \, \sigma ({\cal R}_{CG}^{-1})
= F_{21} L_S^{(-)} F_{12}^{-1}.
$$

\section{Universal {\cal R}-matrix for esoteric quantum algebra}

The ${\cal U}_q(sl(2N+1)$ analogue of the $R$-matrix (\ref{RCG})
\cite{FG1,FG2} has quite a few non-zero entries, so we will not
write it here. This $R$-matrix $R_{FG}$ can be obtained
as the reduction to the fundamental representation of our final result
-- universal ${\cal R}_{FG}$. Instead we start with
the Drinfeld-Jimbo quantum algebra.
The quasitriangular Hopf algebra ${\cal U}_q(sl(2N+1)$ is generated
by $2N$
triples $\{h_i,e_i,f_i\}$ satisfying relations (\ref{triples})
and having coproducts (\ref{coprod}) with $i=1,2,\ldots, 2N$
and three-diagonal Cartan
matrix $\{a_{ij}\}$: $a_{ii}=2$, $a_{ii+1}=a_{i+1i}=-1$.
Serre relations read
$$
e_k^2 e_l - (q+q^{-1})\, e_k e_l e_k + e_l e_k^2,\quad l = k \pm 1,
$$
for two neighbouring root vectors and $e_k e_l = e_l
e_k$
for distant ones. Similar equalities hold for $f_i$ generators.
The ordered product
of $q$-exponentials in the universal $R$-matrix includes factors
corresponding to all positive roots \cite{Ros,KR,LS,KT,T}:
$$
{\cal R}_S = q^{t_0}\prod_{\al\in\Phi^+}^>
\exp_{q^{-2}}(\omega \, e_\al\oo f_\al).
$$
Composite root vectors are defined according to
$$
e_{\al+\bt}=e_\al e_\bt - q e_\bt e_\al, \quad
f_{\al+\bt}=f_\bt f_\al - q^{-1} f_\al f_\bt, \quad \al < \bt \,.
$$
Following the procedure of the preceding section we fix
two Hopf subalgebras ${\cal B}(N)=\{h_j,e_j;j=i=1,2,\ldots,N\}$
and ${\cal B}_-(N)=\{h_k,f_k;k=i=N+1,N+2,\ldots,2N\}$. They
have elements $h_N$ and $h_{N+1}$ commuting nontrivially with
$f_{N+1}$ and $e_N$. Extending ${\cal U}_q(sl(2N+1)$ by the central
element
$C=\sum_{i=1}^{2N+1}e_{ii}$ and performing the "diagonal twist"  with
\be
{\cal F}^{(1)}=\exp(\frac{\gamma}{2}H_{N+1}\wedge Z_{N+1}), \quad
H_{N+1}=\sum_{i=1}^{N+1}e_{ii}\,,  \quad
Z_{N+1}=\sum_{k=N+1}^{2N+1}e_{kk}\,,
\label{Tw1}
\ee
one gets two independent Hopf subalgebras of ${\cal U}_q(gl(2N+1)$.
Their independence can be easily seen from the coproduct on
new generators
$$
e_i\to e_i, \quad i<N,  \quad e_N\to e_N e^{\frac{\gamma}{2}H_{N+1}},
\quad f_{N+1}\to f_{N+1} e^{-\frac{\gamma}{2}Z_{N+1}},\quad
f_i\to f_i, \quad i>N+1,
$$
\be
\Delta_t (e_N)&=&e_N\oo q^{e_{NN}+H_N} +1\oo e_N,
\quad
H_{N}=\sum_{i=1}^{N}e_{ii},
\n\\
\Delta_t (f_{N+1})&=&f_{N+1}\oo 1 +q^{e_{{N+2}{N+2}}+Z_N}\oo f_{N+1},
\quad
Z_{N}=\sum_{k=N+2}^{2N+1}e_{kk}, \n
\ee

Let us introduce "primed" notations $i'=2N+2-i$  and
$\al'_j=\al_{2N+1-j}$, $i,j=1,\ldots N$, corresponding to
reflection of the Dynkin diagram for $sl(2N+1)$.
Now we can identify  the Hopf subalgebra
${\cal B}(N)_-=\{e_{kk},f_{k-1}; k=N+2,\ldots,2N+1\}$
with ${\cal B}(N)_{op}^*$,
the dual to ${\cal B}(N)=\{e_{ii},e_i;i=1,\ldots,N\}$
having the opposite product. The non-vanishing matrix elements
of the pairing between the generators ($\al$ is a simple positive
root)
are
$$
<e_{ii},e_{i'i'}>=1, \quad <e_{\al},f_{\al'}> = 1 \,.
$$
Corresponding canonical element is
given by the ordered product
\be
{\cal F}^{(2)}=\prod_{\al\in\Phi^+}^<
\exp_{q^{2}}(\mu_\al e_\al\oo g_{\al'})
q^{t_0+H_N\oo Z_N},
\label{Tw2}
\ee
where
$\Phi^+$ is the set of all positive roots of $sl(N+1)$ and
$t_0=\sum e_{ii}\oo e_{i'i'}$. Element
$g_{\al'}$ coincides with $f_{\al'}$ if ${\al}$ is a simple root
and is defined by $g_{\al'+\bt'}=g_{\bt'}g_{\al'}-q^{-
1}g_{\al'}g_{\bt'}$
for the case of composite roots. Because of the inverted ordering in
primed roots compared to non-primed ones, $g_{\al'}=f_{\al'}$ only
for the simple
roots. The element ${\cal F}^{(2)}$ involves $N(N+1)/2$
parameters $\mu_\al$
among which those $N$ corresponding to the simple roots are
independent
and $\mu_{\al+\bt}=\mu_{\al}\mu_{\bt}$. Actually, ${\cal F}^{(2)}$
is not exactly the canonical element but related to it by
the transformation of the kind just described (cf. the remark at the
end of Section \ref{Sec2}). We deliberately make no difference between
them so as to simplify the presentation.

Due to the definition of the canonical element it satisfies the
factorization properties (\ref{FP}). Hence one can take it as a twisting
element ${\cal F}^{(2)}$. There is still a part of the centrally extended
Cartan subalgebra $\{e_{ii},i=1,2,\ldots,2N+1\}$
which is invariant with respect to the composition
of two twists ${\cal F}^{(1)}$ and ${\cal F}^{(2)}$.
For $N=1$ its dimension is two:
$
\sum_{i=1}^3 \al_i,([e_{ii},e_{12}]\wedge e_{23}+ e_{12}\wedge
[e_{ii},e_{23}])=0
$
means $2\al_2=\al_1+\al_3$;
in the general case that subalgebra is formed by
$\{e_{ii}-e_{i'i'}, C; i\leq N\}$ ($C$ is the central element).
Hence there is a possibility for additional
diagonal twisting with more parameters in resulting quantum
algebra:
\be
{\cal F}^{(3)}=\exp(A^{ik}(e_{ii}-e_{i'i'})\wedge (e_{kk}-e_{k'k'})
+B^{i}(e_{ii}-e_{i'i'})\wedge C)).
\label{Tw3}
\ee
Finally we arrive to the following

\vspace{0.5in}
\noindent
{\bf Proposition.}
{\em
The esoteric quantum algebra ${\cal U}_{FG}(gl(2N+1))$
defined by the $R$-matrix of the type (\ref{RCG})
\cite{FG1,FG2} is a twisting
of the quasitriangular Hopf algebra ${\cal U}_{q}(gl(2N+1))$
with the twisting element ${\cal F}={\cal F}^{(3)}{\cal F}^{(2)}{\cal
F}^{(1)}$
where ${\cal F}^{(i)}$ are given by expressions (\ref{Tw1}),
(\ref{Tw2}), (\ref{Tw3}), and the universal $R$-matrix
$$
{\cal R}_{FG} = {\cal F}_{21} {\cal R}_{S} {\cal F}^{-1}.
$$
}

\noindent
That way twisted esoteric quantum algebra ${\cal U}_{q}(gl(2N+1))$
and its universal $R$-matrix
have $(N+1)(N+2)/2$ parameters,
which is in accordance with \cite{JC9702028}.

\section{Conclusion}

The explicit expression of the twist ${\cal F}$ has been obtained due to
appropriate choice of initial diagonal twist providing two independent Hopf
subalgebras. We hope that a similar procedure could clarify
interrelation between the Cremmer-Gervais quantum algebra and
${\cal U}_{q}(gl(N))$ for $N > 3$.

There are various possibilities to use the universal/algebraic
twist: \\
(i) relations among the FTR-approach generators of the twisted and
original quantum algebras; \\
(ii) evaluation of the Clebsch-Gordan
coefficients (CGC) of the twisted algebra in terms of the
original CGC and the matrix $F = (\rho_\la \oo \rho_\mu) {\cal F}$
in the tensor product of
the irreducible representations $V_\la \oo V_\mu$; \\
(iii) explicit construction according to the quantum inverse
scattering method of new
integrable models corresponding to twisted $R$-matrices in various
irreducible representations (cf \cite{KS}).

\end{document}